# A simple approach to some Hankel determinants

*Johann Cigler*

Fakultät für Mathematik
Universität Wien
A-1090 Wien, Nordbergstraße 15

johann.cigler@univie.ac.at

**Abstract**
I give simple elementary proofs for some well-known Hankel determinants and their $q-$analogues.

My starting point is the following

**Lemma**
*For given sequences* $s(n)$ *and* $t(n)$ *define* $a(n,k)$ *by*

$$\begin{aligned} &a(0,k)=[k=0] \\ &a(n,0)=s(0)a(n-1,0)+t(0)a(n-1,1) \\ &a(n,k)=a(n-1,k-1)+s(k)a(n-1,k)+t(k)a(n-1,k+1). \end{aligned} \tag{1}$$

*Then the Hankel determinant* $\det\left(a(i+j,0)\right)_{i,j=0}^{n-1}$ *is given by*

$$\det\left(a(i+j,0)\right)_{i,j=0}^{n-1}=\prod_{i=1}^{n-1}\prod_{k=0}^{i-1}t(k). \tag{2}$$

**Proof**
In order to prove this we show first that for all $m,n\in\mathbb{N}$

$$\sum_k a(n,k)a(m,k)\prod_{j=0}^{k-1}t(j)=a(m+n,0). \tag{3}$$

The proof is obvious by induction. For $n=0$ and arbitrary $m$ it is trivially true. If it is true for $n-1$ and arbitrary $m,$ then we get

$$\sum_k a(n,k)a(m,k)\prod_{j=0}^{k-1} t(j) = \sum_k \left(a(n-1,k-1)+s(k)a(n-1,k)+t(k)a(n-1,k+1)\right)a(m,k)\prod_{j=0}^{k-1} t(j)$$
$$= \sum_k a(n-1,k)\left(a(m,k+1)t(k)+s(k)a(m,k)+a(m,k-1)\right)\prod_{j=0}^{k-1} t(j)$$
$$= \sum_k a(n-1,k)a(m+1,k)\prod_{j=0}^{k-1} t(j) = a(n+m,0).$$

Consider now the matrices

$A_n = \left(a(i,j)\right)_{i,j=0}^{n-1}$,

$D_n = \left([i=j]\prod_{k=0}^{i-1} t(j)\right)_{i,j=0}^{n-1}$,

and

$H_n = \left((a(i+j,0)\right)_{i,j=0}^{n-1}$.

Then (3) implies

$$A_n D_n A_n^t = H_n.$$

Therefore we get

$$\det\left(a(i+j,0)\right)_{i,j=0}^{n-1} = \det H_n = \prod_{i=1}^{n-1}\prod_{k=0}^{i-1} t(k). \qquad (4)$$

**Remark**

This Lemma is well- known and intimately connected with the approach to Hankel determinants via orthogonal polynomials (see [6]). It is especially useful if for a given sequence $\left(a(n,0)\right)_{n\geq 0}$ explicit expressions for $s(n), t(n)$ and $a(n,k)$ exist. To find such expressions it is often convenient to compute the first values of the orthogonal polynomials $p(n,x)$ (cf. e.g. [4] (1.10)) and their Favard resolution [4] (1.11) and try to guess $s(n)$ and $t(n)$ and the explicit form of $a(n,k)$. Other approaches to this lemma can be found in [1] and [8].

If all $s(n)=0$ identity (1) reduces to

$$\begin{aligned} &A(0,k) = [k=0] \\ &A(n,0) = T(0)A(n-1,1) \\ &A(n,k) = A(n-1,k-1) + T(k)A(n-1,k+1). \end{aligned} \qquad (5)$$

In this case $A(2n,2k+1) = A(2n+1,2k) = 0$ for all $n,k$.

If we define

$$a(n,k) = A(2n,2k), \qquad (6)$$

then it is easily verified that (1) holds with

$$\begin{aligned} s(0) &= T(0), \\ s(n) &= T(2n-1) + T(2n), \\ t(n) &= T(2n)T(2n+1). \end{aligned} \qquad (7)$$

Therefore it is convenient to consider first sequences of the form $\left(c(0), 0, c(1), 0, c(2), 0, \cdots\right)$.

Let us first look at the situation for the famous example of the Catalan numbers

$$c(n) = C_n = \frac{1}{n+1}\binom{2n}{n}. \qquad (8)$$

The table $\left(A(n,k)\right)$ can be found in [5], A053121. The first terms are

| | | | | | | | |
|---|---|---|---|---|---|---|---|
| 1 | | | | | | | |
| 0 | 1 | | | | | | |
| 1 | 0 | 1 | | | | | |
| 0 | 2 | 0 | 1 | | | | |
| 2 | 0 | 3 | 0 | 1 | | | |
| 0 | 5 | 0 | 4 | 0 | 1 | | |
| 5 | 0 | 9 | 0 | 5 | 0 | 1 | |
| 0 | 14 | 0 | 14 | 0 | 6 | 0 | 1 |

It is well known that $T(n) = 1$ and $A(n,k) = \frac{k+1}{n+1}\binom{n+1}{\frac{n-k}{2}}$, if $n \equiv k \bmod 2$ and $0$ else. This can also be written in the form $A(2n,2k) = \binom{2n}{n-k} - \binom{2n}{n-k-1} = \frac{2k+1}{n+k+1}\binom{2n}{n-k}$ and $A(2n+1,2k+1) = \binom{2n+1}{n-k} - \binom{2n+1}{n-k-1} = \frac{2k+2}{n+k+2}\binom{2n+1}{n-k}$.

To prove this fact it suffices to verify (5). This is almost trivial:

$A(n,0) = T(0)A(n-1,1)$ is equivalent with

$$A(2n,0) = \binom{2n}{n} - \binom{2n}{n-1} = A(2n-1,1) = \binom{2n-1}{n-1} - \binom{2n-1}{n-2}$$

or

$$\binom{2n}{n}-\binom{2n-1}{n-1}=\binom{2n}{n-1}-\binom{2n-1}{n-2}=\binom{2n-1}{n}.$$

$A(n,k)=A(n-1,k-1)+T(k)A(n-1,k+1)$
reduces to
$A(2n,2k)=A(2n-1,2k-1)+A(2n-1,2k+1)$
and
$A(2n+1,2k+1)=A(2n,2k)+A(2n,2k+2).$
The first identity is

$$\binom{2n}{n-k}-\binom{2n}{n-k-1}=\binom{2n-1}{n-k}-\binom{2n-1}{n-k-1}+\binom{2n-1}{n-k-1}-\binom{2n-1}{n-k-2}$$

and the second one

$$\binom{2n+1}{n-k}-\binom{2n+1}{n-k-1}=\binom{2n}{n-k}-\binom{2n}{n-k-1}+\binom{2n}{n-k-1}-\binom{2n}{n-k-2}.$$

Both follow from the recurrences of the binomial coefficients.

For the sequence $a(n,0)=C_n$ we get $a(n,k)=\binom{2n}{n-k}-\binom{2n}{n-k-1}=\frac{2k+1}{n+k+1}\binom{2n}{n-k}$,
$s(0)=1, s(n)=2$ for $n>0$ and $t(n)=1.$
This gives the table [5], A039599

| | | | | |
|---|---|---|---|---|
| 1 | | | | |
| 1 | 1 | | | |
| 2 | 3 | 1 | | |
| 5 | 9 | 5 | 1 | |
| 14 | 28 | 20 | 7 | 1 |

An immediate consequence of (4) is the well-known result $\det\left(C_{i+j}\right)_{i,j=0}^{n-1}=1$ for all $n\in\mathbb{N}.$

It is also easy to verify that

$$\sum_{k=0}^{n}(-1)^k a(n,k)=[n=0] \tag{9}$$

and

$$\sum_{k=0}^{n}a(n,k)=\binom{2n}{n}. \tag{10}$$

But we cannot immediately deduce that $\det\left(C_{i+j+1}\right)_{i,j=0}^{n-1}=1$ or obtain the higher Hankel determinants $\det\left(C_{i+j+m}\right)_{i,j=0}^{n-1}$.

But observe that $C_n = 4^n \prod_{j=0}^{n-1} \frac{1+2j}{4+2j}$ and therefore $C_{n+m} = 4^{n+m} \prod_{j=0}^{m-1} \frac{1+2j}{4+2j} \prod_{j=0}^{n-1} \frac{1+2m+2j}{4+2m+2j}.$

If we can find explicit values of $s(n), t(n)$ and $a(n,k)$ for $c(n) = \prod_{j=0}^{n-1} \frac{1+2m+2j}{4+2m+2j}$, then we can compute these higher Hankel determinants. This is indeed the case, even for the more general case $c(n) = \frac{\prod_{j=0}^{n-1}(b+jc)}{\prod_{j=0}^{n-1}(a+jc)}$. It also remains true for $c(n) = \prod_{j=0}^{n-1} \frac{1-q^{b+jc}}{1-q^{a+jc}}$, which for $q \to 1$ tends to $\frac{\prod_{j=0}^{n-1}(b+jc)}{\prod_{j=0}^{n-1}(a+jc)}$ and more generally for $c(n) = \prod_{j=0}^{n-1} \frac{1-bq^j}{1-aq^j}$ for arbitrary $a, b.$

We shall later use the trivial fact that (1) is equivalent with

$$
\begin{aligned}
&a(0,k) = [k=0] \\
&\left(x^n a(n,0)\right) = \left(xs(0)\right)\left(x^{n-1}a(n-1,0)\right) + \left(xt(0)\right)\left(x^{n-1}a(n-1,1)\right) \\
&\left(x^{n-k} a(n,k)\right) = \left(x^{n-k}a(n-1,k-1)\right) + \left(xs(k)\right)\left(x^{n-1-k}a(n-1,k)\right) + \left(x^2 t(k)\right)\left(x^{n-k-2}a(n-1,k+1)\right).
\end{aligned}
\tag{11}
$$

I will not show here how to guess the following results (this the reader can do for himself using the above hint), but will only verify the relevant identities.

We shall use the usual $q-$notations, e.g. $(x;q)_n = \prod_{j=0}^{n-1}(1-q^j x)$, $[n] = \frac{1-q^n}{1-q}$ and $\begin{bmatrix} n \\ k \end{bmatrix} = \frac{(q;q)_n}{(q;q)_k (q;q)_{n-k}}$ for $0 \le k \le n$ and $\begin{bmatrix} n \\ k \end{bmatrix} = 0$ for $k < 0$ and $k > n.$

**Theorem 1**

*Let*

$$c(n,a,b,q) = \frac{(b;q)_n}{(a;q)_n}. \tag{12}$$

*Condition (5) is satisfied if we set*

$$
\begin{aligned}
T(2n,a,b,q) &= \frac{q^n(1-q^n b)(1-q^{n-1}a)}{(1-q^{2n-1}a)(1-q^{2n}a)}, \\
T(2n+1,a,b,q) &= \frac{q^n(1-q^{n+1})(b-q^n a)}{(1-q^{2n+1}a)(1-q^{2n}a)}
\end{aligned}
\tag{13}
$$

*and*

$$A(2n,2k,a,b,q)=\begin{bmatrix} n \\ k \end{bmatrix} c(n-k,q^{2k}a,q^{k}b,q),$$
$$A(2n+1,2k+1,a,b,q)=\begin{bmatrix} n \\ k \end{bmatrix} c(n-k,q^{2k+1}a,q^{k+1}b,q). \tag{14}$$

**Proof**

We have to show that the identities

$$A(2n+2,2k,a,b,q)-A(2n+1,2k-1,a,b,q)-T(2k,a,b,q)A(2n+1,2k+1,a,b,q)=0 \tag{15}$$

and

$$A(2n+1,2k+1,a,b,q)-A(2n,2k,a,b,q)-T(2k+1,a,b,q)A(2n,2k+2,a,b,q)=0 \tag{16}$$

hold.
The left-hand side of (15) is

$$\begin{bmatrix} n+1 \\ k \end{bmatrix} c(n-k+1,q^{2k}a,q^{k}b,q)-\begin{bmatrix} n \\ k-1 \end{bmatrix} c(n-k+1,q^{2k-1}a,q^{k}b,q)$$
$$-T(2k,a,b,q)\begin{bmatrix} n \\ k \end{bmatrix} c(n-k,q^{2k+1}a,q^{k+1}b,q)$$
$$=\frac{(q;q)_{n+1}}{(q;q)_k\,(q;q)_{n+1-k}}\frac{\left(q^kb;q\right)_{n-k+1}}{\left(q^{2k}a;q\right)_{n-k+1}}-\frac{(q;q)_n}{(q;q)_{k-1}\,(q;q)_{n+1-k}}\frac{\left(q^kb;q\right)_{n-k+1}}{\left(q^{2k-1}a;q\right)_{n-k+1}}$$
$$-\frac{q^k(1-q^{k+1})(1-q^{k-1}a)}{(1-q^{2k-1}a)(1-q^{2k}a)}\frac{(q;q)_n}{(q;q)_k\,(q;q)_{n-k}}\frac{\left(q^{k+1}b;q\right)_{n-k}}{\left(q^{2k+1}a;q\right)_{n-k}}$$
$$=\frac{(q;q)_n}{(q;q)_k\,(q;q)_{n+1-k}}\frac{\left(q^kb;q\right)_{n-k+1}}{\left(q^{2k-1}a;q\right)_{n-k+2}}\left(\left(1-q^{n+1}\right)\left(1-q^{2k-1}a\right)-\left(1-q^k\right)\left(1-q^{n+k}a\right)-q^k\left(1-q^{k-1}a\right)\left(1-q^{n+1-k}\right)\right)$$

Since

$$\left(1-q^{n+1}\right)\left(1-q^{2k-1}a\right)-\left(1-q^k\right)\left(1-q^{n+k}a\right)-q^k\left(1-q^{k-1}a\right)\left(1-q^{n+1-k}\right)=0,$$

we see that (15) is true.

In the same way the left-hand side of (16) reduces to

$$\begin{bmatrix} n \\ k \end{bmatrix} c(n-k, q^{2k+1}a, q^{k+1}b, q) - \begin{bmatrix} n \\ k \end{bmatrix} c(n-k, q^{2k}a, q^k b, q)$$

$$-\frac{q^k(1-q^{k+1})(b-q^k a)}{(1-q^{2k+1}a)(1-q^{2k}a)} \begin{bmatrix} n \\ k+1 \end{bmatrix} c(n-k-1, q^{2k+2}a, q^{k+1}b, q)$$

$$= \frac{(q;q)_n}{(q;q)_k (q;q)_{n-k}} \frac{(q^{k+1}b;q)_{n-k}}{(q^{2k+1}a;q)_{n-k}} - \frac{(q;q)_n}{(q;q)_k (q;q)_{n-k}} \frac{(q^k b;q)_{n-k}}{(q^{2k}a;q)_{n-k}}$$

$$-\frac{q^k(1-q^{k+1})(b-q^k a)}{(1-q^{2k-1}a)(1-q^{2k}a)} \frac{(q;q)_n}{(q;q)_{k+1}(q;q)_{n-k-1}} \frac{(q^{k+1}b;q)_{n-k-1}}{(q^{2k+2}a;q)_{n-k-1}}$$

$$= \frac{(q;q)_n}{(q;q)_k (q;q)_{n-k}} \frac{(q^{k+1}b;q)_{n-k-1}}{(q^{2k}a;q)_{n-k+1}} \Big( (1-q^n b)(1-q^{2k}a) - (1-q^k b)(1-q^{n+k}a) - q^k (b-q^k a)(1-q^{n-k}) \Big) = 0.$$

This implies

**Theorem 2**
*Let*

$$c(n,a,b,q) = \frac{(b;q)_n}{(a;q)_n}. \tag{17}$$

*Then the Hankel determinants*

$$d(n,m,a,b,q) = \det\big(c(i+j+m,a,b,q)\big)_{i,j=0}^{n-1} \tag{18}$$

*are given by*

$$d(n,0,a,b,q) = q^{2\binom{n}{3}} \prod_{k=1}^{n-1} \frac{(b;q)_k (q;q)_k \prod_{j=0}^{k-1} (b-q^j a)}{(q^{k-1}a;q)_k (a;q)_{2k}} \tag{19}$$

*and*

$$d(n,m,a,b,q) = d(n,0,a,b,q) q^{m\binom{n}{2}} \prod_{j=0}^{m-1} \frac{(q^j b;q)_n}{(q^{n-1+j}a;q)_n}. \tag{20}$$

**Proof**

By (2) we have $d(n,0,a,b,q) = \prod_{k=1}^{n-1} \prod_{j=0}^{k-1} t(j,a,b,q) = \prod_{k=1}^{n-1} \prod_{j=0}^{k-1} T(2j,a,b,q) T(2j+1,a,b,q).$

From

$$\prod_{j=0}^{k-1} T(2j,a,b,q) = q^{\binom{k}{2}} \frac{(b;q)_k}{(q^{k-1}a;q)_k}$$

and

$$\prod_{j=0}^{k-1} T(2j+1,a,b,q) = q^{\binom{k}{2}} \frac{(q;q)_k}{(a;q)_{2k}} \prod_{j=0}^{k-1} (b-q^j a)$$

we deduce (19).

By (19) we get

$$\frac{d(n,0,qa,qb,q)}{d(n,0,a,b,q)} = q^{\binom{n}{2}} \left(\frac{1-a}{1-b}\right)^n \frac{(b;q)_n}{(q^{n-1}a;q)_n}.$$

This implies

$$d(n,1,a,b,q) = \det\left(\frac{(b;q)_{i+j+1}}{(a;q)_{i+j+1}}\right)_{i,j=0}^{n-1} = \det\left(\frac{(1-b)}{(1-a)}\frac{(qb;q)_{i+j}}{(qa;q)_{i+j}}\right) = \left(\frac{1-b}{1-a}\right)^n d(n,0,qa,qb,q)$$

$$= q^{\binom{n}{2}} \frac{(b;q)_n}{(q^{n-1}a;q)_n} d(n,0,a,b,q).$$

Iterating this argument we get

$$d(n,m,a,b,q) = \left(\frac{(b;q)_m}{(a;q)_m}\right)^n d(n,0,q^m a,q^m b,q)$$

$$= q^{m\binom{n}{2}} \prod_{j=0}^{m-1} \frac{(q^j b;q)_n}{(q^{n-1+j}a;q)_n} d(n,0,a,b,q).$$

If we change $a \to q^a, b \to q^b, q \to q^c$ and let $q \to 1$ we get

**Corollary 1**

*Let*

$$u(n,a,b,c) = \frac{\prod_{j=0}^{n-1}(b+jc)}{\prod_{j=0}^{n-1}(a+jc)}. \tag{21}$$

*Condition (5) is satisfied if we set*

$$\begin{aligned} T(2n,a,b,c) &= \frac{(a+(n-1)c)(b+nc)}{(a+2nc)(a+(2n-1)c)}, \\ T(2n+1,a,b,c) &= \frac{(n+1)c(a-b+nc)}{(a+2nc)(a+(2n+1)c)} \end{aligned} \tag{22}$$

*and*

$$\begin{aligned} A(2n,2k,a,b,c) &= \binom{n}{k} u(n-k,a+2kc,b+kc,c), \\ A(2n+1,2k+1,a,b,c) &= \binom{n}{k} u(n-k,a+(2k+1)c,b+(k+1)c,c). \end{aligned} \tag{23}$$

**Corollary 2**

*Let*

$$u(n,a,b,c)=\frac{\prod_{j=0}^{n-1}\left(b+jc\right)}{\prod_{j=0}^{n-1}\left(a+jc\right)}. \tag{24}$$

*Then the Hankel determinants*

$$D(n,m,a,b,c)=\det\left(u(i+j+m,a,b,c)\right)_{i,j=0}^{n-1} \tag{25}$$

*are given by*

$$D(n,0,a,b,c)=\prod_{k=1}^{n-1}\frac{k!c^k\prod_{j=0}^{k-1}(b+jc)(a-b+jc)}{\prod_{j=0}^{k-1}\left(a+(j+k-1)c\right)\prod_{j=0}^{2k-1}\left(a+jc\right)} \tag{26}$$

*and*

$$D(n,m,a,b,c)=D(n,0,a,b,c)\prod_{j=0}^{m-1}\prod_{i=0}^{n-1}\frac{\left(b+(j+i)c\right)}{\left(a+(i+n+j-1)c\right)}. \tag{27}$$

These Corollaries can of course also be proved directly.

## Hankel determinants of Catalan numbers

Now let us consider again the Catalan numbers. Instead of $C_n$ we study

$$u(n,4,1,2)=\prod_{j=0}^{n-1}\frac{1+2j}{4+2j}. \tag{28}$$

Here we get $T(n)=\frac{1}{4}$ and find again that

$$A(2n,2k)=\binom{n}{k}u(n-k,4+4k,1+2k,2)=\binom{n}{k}\frac{(1+2k)(3+2k)\cdots(2n-1)}{(4+4k)(6+4k)\cdots(2n+2k+2)}$$
$$=\frac{1}{4^{n-k}}\frac{n!}{k!(n-k)!}\frac{(2n)!}{(2k)!(2+2k)(4+2k)\cdots(2n)}\frac{(2k+1)!}{(n+k+1)!}=\frac{1}{4^{n-k}}\frac{2k+1}{n+k+1}\binom{2n}{n-k}$$

and

$$A(2n+1,2k+1)=\frac{1}{4^{n-k}}\frac{2k+2}{n+k+2}\binom{2n+1}{n-k}.$$

It is evident that $\det\left(u(i+j,4,1,2)\right)=\frac{1}{4^{\binom{n}{2}}}$.

In order to compute the higher Hankel determinants consider the factor
$\prod_{i=0}^{n-1} \frac{(b+(j+i)c)}{(a+(i+n+j-1)c)}$ which occurs in (27).

This reduces in our case to

$$\prod_{i=0}^{n-1} \frac{1+2j+2i}{4+2i+2n+2j-2} = \frac{1}{2^n} \prod_{i=0}^{n-1} \frac{1+2j+2i}{1+n+j+i} = \frac{1}{2^n} \frac{(2j+2n)!}{(2j)! \prod_{\ell=1}^{n} (2j+2\ell)} \frac{(n+j)!}{(2n+j)!}$$

$$= \frac{1}{4^n} \frac{(2j+2n)!\, j!}{(2j)!(2n+j)!} = \frac{1}{4^n} \prod_{i=1}^{j} \frac{2n+j+i}{j+i}.$$

Thus we get

$$D(n,m,4,1,2) = \frac{1}{4^{\binom{n}{2}+mn}} \prod_{j=0}^{m-1} \prod_{i=1}^{j} \frac{2n+j+i}{j+i}.$$

Therefore we get for the Catalan numbers themselves the well-known result (cf. [7], Theorem 33 and the historical comments given there)

$$\det\left(C_{i+j+m}\right)_{i,j=0}^{n-1} = \prod_{j=1}^{m-1} \prod_{i=1}^{j} \frac{2n+j+i}{j+i}. \tag{29}$$

Since $\det\left(C_{i+j+1}\right)_{i,j=0}^{n-1} = 1,$ we get also the well-known result, that the Catalan numbers are the uniquely determined numbers satisfying $\det\left(C_{i+j}\right)_{i,j=0}^{n-1} = 1$ and $\det\left(C_{i+j+1}\right)_{i,j=0}^{n-1} = 1.$

**Some applications of Theorem 2**

1) Let $c(n) = c(n,0,x,q) = (x;q)_n$.

Then we get

$$\det\left((x;q)_{i+j+m}\right)_{i,j=0}^{n-1} = q^{2\binom{n}{3}+m\binom{n}{2}} x^{\binom{n}{2}} \prod_{k=0}^{n-1} (x;q)_{k+m} (q;q)_k. \tag{30}$$

As special case we note

$$\det\left([i+j+m]!\right)_{i,j=0}^{n-1} = q^{2\binom{n}{3}+(m+1)\binom{n}{2}} \prod_{k=0}^{n-1} [k+m]![k]!. \tag{31}$$

If we define

$$\langle x \rangle_n := \prod_{j=0}^{n-1} \left(x-[j]\right), \tag{32}$$

then we get $\langle x\rangle_n = \frac{(1+(q-1)x)^n}{(q-1)^n}\left(\frac{1}{1+(q-1)x};q\right)_n$ .

Therefore we have

$$\det\left(\langle x\rangle_{i+j+m}\right)_{i,j=0}^{n-1} = (-1)^{\binom{n}{2}} q^{2\binom{n}{3}+m\binom{n}{2}} \prod_{j=0}^{n-1} [j]!\langle x\rangle_{j+m}. \quad (33)$$

2) Let $c(n) = c(n,q,q^{d+1},q) = \frac{\left(q^{d+1};q\right)_n}{\left(q;q\right)_n}$.

Here we get

$$d(n,0,q,q^{d+1},q) = q^{2\binom{n}{3}} \prod_{k=1}^{n-1} \frac{\left(q^{d+1};q\right)_k \left(q;q\right)_k \prod_{j=0}^{k-1}\left(q^{d+1}-q^{j+1}\right)}{\left(q^k;q\right)_k \left(q;q\right)_{2k}}$$

$$= q^{2\binom{n}{3}+\binom{n}{2}} \prod_{k=1}^{n-1} \frac{\left(q^{d+1};q\right)_k \prod_{j=0}^{k-1}\left(q^d-q^j\right)}{\left(q^k;q\right)_k \left(q^{k+1};q\right)_k}.$$

For $d = m \in \mathbb{N}$ we get $c(n) = \begin{bmatrix} n+m \\ m \end{bmatrix}$.

The Hankel determinant reduces to

$$\det\left(\begin{bmatrix} i+j+m \\ m \end{bmatrix}\right)_{i,j=0}^{n-1} = (-1)^{\binom{n}{2}} q^{\frac{n(n-1)^2}{2}} \prod_{k=0}^{n-1} \frac{\begin{bmatrix} m+k \\ 2k \end{bmatrix}}{\begin{bmatrix} 2k-1 \\ k \end{bmatrix}}. \quad (34)$$

This result is due to Carlitz [3]. For some generalizations see Krattenthaler [6], especially Theorem 26.

3) Choose $c(n) = c(n,2,1,1) = \frac{1-q}{1-q^{n+1}} = \frac{1}{[n+1]}$.

Then (20) gives a $q-$analogue of the determinant of the Hilbert matrix

$$\det\left(\frac{1}{[i+j+m+1]}\right)_{i,j=0}^{n-1} = q^{m\binom{n}{2}+\frac{n(n-1)(2n-1)}{6}} \prod_{j=0}^{m-1} \frac{\left([j+n]!\right)^2}{[j]![2n+j]!} \prod_{j=0}^{n-1} \frac{\left([j]!\right)^3}{[n+j]!}. \quad (35)$$

4) Let $c(n)=(1-q)^n c(n,1,0,1)=\frac{(1-q)^n}{(q;q)_n}=\frac{1}{[n]!}$.

Then

$$\det\left(\frac{1}{[i+j+m]!}\right)_{i,j=0}^{n-1}=(-1)^{\binom{n}{2}}q^{\frac{n(n-1)^2}{2}}\prod_{j=0}^{n+m-2}\frac{[j]!}{[n+j]!}. \tag{36}$$

For $\det\left(\frac{1}{[i+j]!}\right)_{i,j=0}^{n-1}=(-1)^{\binom{n}{2}}q^{\frac{n(n-1)^2}{2}}\prod_{j=0}^{n-2}\frac{[j]!}{[n+j]!}$

and

$$\det\left(\frac{1}{[i+j+m]!}\right)_{i,j=0}^{n-1}=q^{m\binom{n}{2}}\prod_{j=0}^{m-1}\frac{[n+j-1]!}{[2n+j-1]!}\det\left(\frac{1}{[i+j]!}\right)_{i,j=0}^{n-1}$$
$$=(-1)^{\binom{n}{2}}q^{\frac{n(n-1)^2}{2}}\prod_{j=0}^{n+m-2}\frac{[j]!}{[n+j]!}.$$

5) Now we choose

$$c(n)=c(n,q^2,q,q^2)=\prod_{j=1}^{n}\frac{1-q^{2j-1}}{1-q^{2j}}=\begin{bmatrix}2n\\n\end{bmatrix}\frac{1}{\prod_{j=1}^{n}(1+q^j)^2}=(-1)^n q^{n^2}\begin{bmatrix}-\frac{1}{2}\\n\end{bmatrix}_{q^2}.$$

Here we get

$$d(n,0,q^2,q,q^2)=q^{\frac{n(n-1)(4n-5)}{6}}\frac{1}{\prod_{j=1}^{2n-2}(1+q^j)^{2n-1-j}} \tag{37}$$

and

$$d(n,m,q^2,q,q^2)=q^{2m\binom{n}{2}}\prod_{j=0}^{m-1}\frac{\left(q^{2j+1};q^2\right)_n}{\left(q^{2n+2j};q^2\right)_n}d(n,0,q^2,q,q^2). \tag{38}$$

The right-hand side can be simplified:

$$\prod_{j=0}^{m-1}\frac{\left(q^{2j+1};q^2\right)_n}{\left(q^{2n+2j};q^2\right)_n}=\prod_{j=0}^{m-1}\prod_{i=0}^{n-1}\frac{[2j+2i+1]}{[2j+2i+2n]}=\prod_{j=0}^{m-1}\frac{[2j+2n-1]!}{[2][4]\cdots[4n+2j-2][1][3]\cdots[2j-1]}$$
$$=\prod_{j=0}^{m-1}\frac{[2n+2j-1]!}{[2n+j-1]![1][3]\cdots[2j-1](1+q)(1+q^2)\cdots(1+q^{2n+j-1})}$$
$$=\prod_{j=0}^{m-1}\frac{[2n+j]\cdots[2n+2j-1]}{[1][3]\cdots[2j-1](1+q)(1+q^2)\cdots(1+q^{2n+j-1})}.$$

Observe that

$$\frac{[j+1][j+2]\ldots[2j]}{[1][3]\cdots[2j-1]}=\frac{[2][4]\cdots[2j][j+1][j+2]\ldots[2j]}{[2j]!}=(1+q)\cdots(1+q^{j}).$$

This implies

$$\prod_{j=0}^{m-1}\frac{[2n+j]\cdots[2n+2j-1]}{[1][3]\cdots[2j-1](1+q)(1+q^2)\cdots(1+q^{2n+j-1})}=\prod_{j=0}^{m-1}\frac{[2n+j]\cdots[2n+2j-1](1+q)(1+q^2)\cdots(1+q^{j})}{(1+q)(1+q^2)\cdots(1+q^{2n+j-1})[j+1]\cdots[2j]}$$
$$=\prod_{j=0}^{m-1}\frac{[2n+j]\cdots[2n+2j-1]}{(1+q^{j+1})(1+q^{j+2})\cdots(1+q^{2n+j-1})[j+1]\cdots[2j]}=\prod_{j=0}^{m-1}\frac{1}{\left(-q^{j+1};q\right)_{2n-1}}\prod_{i=1}^{j}\frac{[2n+j+i-1]}{[j+i]}.$$

Finally we get

$$d(n,m,q^2,q,q^2)=q^{2m\binom{n}{2}}\prod_{j=0}^{m-1}\frac{1}{\left(-q^{j+1};q\right)_{2n-1}}\prod_{i=1}^{j}\frac{[2n+j+i-1]}{[j+i]}d(n,0,q^2,q,q^2). \qquad (39)$$

We note also that

$$T(0,q^2,q,q^2)=\frac{1}{1+q} \qquad (40)$$

and

$$T(n,q^2,q,q^2)=\frac{q^{n}}{(1+q^{n})(1+q^{n+1})} \qquad (41)$$

for $n>0$.

If we let $q\to 1$ we obtain $c(n)\to\frac{1}{4^{n}}\binom{2n}{n}$. This implies that the Hankel determinants of the central binomial coefficients are given by

$$\det\left(\binom{2i+2j+2m}{i+j+m}\right)_{i,j=0}^{n-1}=2^{n-1+m}\prod_{j=0}^{m-1}\prod_{i=1}^{j}\frac{2n+j+i-1}{j+i}. \qquad (42)$$

Similar determinants have also been obtained with other methods by Krattenthaler in [6], especially Theorem 26.

From (40) and (41) and (7) we get for $c(n)$ the values $s(n)=\frac{1}{2}$, $t(0)=\frac{1}{8}$ and $t(n)=\frac{1}{4^2}$ for $n>0.$

Therefore we get for $a(n,0)=\binom{2n}{n}=4^{n}c(n)$ using (11) that $s(n)=4\cdot\frac{1}{2}=2$, $t(0)=4^2\cdot\frac{1}{8}=2$ and $t(n)=4^2\cdot\frac{1}{4^2}=1$ for $n>0.$

This leads to the table [5], A094527,

| | | | | |
|---|---|---|---|---|
| 1 | | | | |
| 2 | 1 | | | |
| 6 | 4 | 1 | | |
| 20 | 15 | 6 | 1 | |
| 70 | 56 | 28 | 8 | 1 |

Here it is easily verified that $a(n,k)=\binom{2n}{n-k}$.

Observe that

$$(1+q)c(n+1,q^2,q,q^2)=c(n,q^4,q^3,q^2)=\frac{(q^3;q^2)_n}{(q^4;q^2)_n}=\begin{bmatrix}2n+1\\ n\end{bmatrix}\frac{1+q}{(1+q^{n+1})}\frac{1}{\prod_{j=1}^{n}(1+q^j)^2}.$$

Therefore

$$\det\left(c(i+j+m,q^4,q^3,q^2)\right)=(1+q)^{\binom{n}{2}}\det\left(c(i+j+m+1,q^2,q,q^2)\right).$$

This implies that the Hankel determinants for the binomial coefficients $\binom{2n+1}{n}$ are

$$\det\left(\binom{2i+2j+2m+1}{i+j+m}\right)_{i,j=0}^{n-1}=\frac{1}{2^n}\det\left(\binom{2i+2j+2m+2}{i+j+m+1}\right)_{i,j=0}^{n-1}. \tag{43}$$

6) Another interesting case is the following $q-$analogue of the Catalan numbers which has first been studied by George Andrews [2]:

$$c(n)=c(n,q^4,q,q^2)=(-1)^n q^{n^2}(1+q)\begin{bmatrix}\frac{1}{2}\\ n+1\end{bmatrix}_{q^2}=\frac{1}{[n+1]}\begin{bmatrix}2n\\ n\end{bmatrix}\frac{1+q}{(1+q^{n+1})\prod_{j=1}^{n}(1+q^j)^2}.$$

Here we get

$$d(n,0,q^4,q,q^2)=\frac{q^{\frac{n(n-1)(4n-5)}{6}}}{(1+q)^{n-1}\prod_{j=0}^{2n-3}(1+q^{j+2})^{2n-2-j}} \tag{44}$$

and with a similar argument as above

$$d(n,m,q^4,q,q^2)=q^{2m\binom{n}{2}}\prod_{j=0}^{m-1}\frac{1}{\left(-q^{j+1};q\right)_{2n}}\prod_{i=1}^{j}\frac{[2n+j+i]}{[j+i]}d(n,0,q^4,q,q^2). \tag{45}$$

Letting $q\to 1$ and simplifying we obtain again (29).

**Some identities connected with this method**

As a generalization of (9) we show, that with $A(n,k)$ defined by (5) the following identity holds:

$$\sum_{k=0}^{n}(-1)^k A(2n,2k)\prod_{j=0}^{k-1}T(2j)=[n=0]. \tag{46}$$

For the proof let $T(-1)=0$ and observe that $A(2n,-2)=A(2n,2n+2)=0$. Then we get

$$\sum_{k=0}^{n+1}(-1)^k A(2n+2,2k)\prod_{j=0}^{k-1}T(2j)$$

$$=\sum_{k=0}^{n+1}(-1)^k \left(A(2n,2k-2)+\left(T(2k-1)+T(2k)\right)A(2n,2k)+T(2k)T(2k+1)a(2n,2k+2)\prod_{j=0}^{k-1}T(2j)\right.$$

$$=\sum_{k=0}^{n}(-1)^k A(2n,2k)\prod_{j=0}^{k-1}T(2j)\left(T(2k-1)+T(2k)-T(2k)-T(2k-1)\right)=0.$$

If we choose $c(n,a,b,q)=\dfrac{(b;q)_n}{(a;q)_n}$ then formula (46) is equivalent to

$$\sum_{k=0}^{n}(-1)^k q^{\binom{k}{2}}\begin{bmatrix}n\\k\end{bmatrix}\left(1-q^{2k}a\right)\frac{1}{\left(q^k a;q\right)_{n+1}}=[n=0]. \tag{47}$$

For the left-hand side of (46) is

$$\sum_{k=0}^{n}(-1)^k\begin{bmatrix}n\\k\end{bmatrix}\frac{\left((q^k b;q)_{n-k}\right)}{\left((q^{2k}a;q)_{n-k}\right)}q^{\binom{k}{2}}\frac{(b;q)_k)}{\left((q^{k-1}a;q)_k\right)}$$

$$=\sum_{k=0}^{n}(-1)^k q^{\binom{k}{2}}\begin{bmatrix}n\\k\end{bmatrix}(1-q^{2k-1}a)\frac{\left((b;q)_n\right)}{\left((q^{k-1}a;q)_{n+1}\right)}$$

Now replace $a\to qa$ and factor through $(b;q)_n$.

A companion to (47) is

$$\sum_{k=0}^{n}q^{\binom{k}{2}}\begin{bmatrix}n\\k\end{bmatrix}\frac{(1-q^{2k}a)}{\left(q^k a;q\right)_{n+1}}=\frac{\left(-1;q\right)_n}{\left(qa;q^2\right)_n}. \tag{48}$$

This is easily proved using Zeilberger's algorithm. E.g. qZeil gives

```
qZeil[ q^Binomial[k, 2] qBinomial[n, k, q] (1 - a q^(2 k))
  qPochhammer[q a, q^2, n] / qPochhammer[q^k a, q, n + 1] , {k, 0, n}, n, 1]
```

```
SUM[n] == (1 + q^(-1+n)) SUM[-1 + n]
```

The proof depends on the following certificate

```
Cert[]
```

$$\frac{q^{-k+n}\left(-1+a\,q^{k}\right)\left(-1+a\,q^{1+2k}\right)\left(-q^{k}+q^{n}\right)}{\left(-1+a\,q^{2k}\right)\left(-1+q^{n}\right)\left(q-a\,q^{2n}\right)}$$

Identity (48) is equivalent with

$$\sum_{k=0}^{n} A(2n,2k,a,b,q)\prod_{j=0}^{k-1} T(2j,a,b,q)=\frac{(b;q)_n(-1;q)_n}{(a;q^2)_n}. \qquad (49)$$

For the left-hand side is

$$\sum_{k=0}^{n}\begin{bmatrix} n\\ k\end{bmatrix} c(n-k,q^{2k}a,q^k b,q)q^{\binom{k}{2}}\frac{(b;q)_k}{(q^{k-1}a;q)_k}=\sum_{k=0}^{n}\begin{bmatrix} n\\ k\end{bmatrix}\frac{(q^k b;q)_{n-k}}{(q^{2k}a;q)_{n-k}}q^{\binom{k}{2}}\frac{(b;q)_k}{(q^{k-1}a;q)_k}$$

$$=\sum_{k=0}^{n}\begin{bmatrix} n\\ k\end{bmatrix} q^{\binom{k}{2}}\frac{(b;q)_n(1-q^{2k-1}a)}{(q^{k-1}a;q)_{n+1}}=\frac{(b;q)_n(-1;q)_n}{(a;q^2)_n}.$$

An interesting special case is the following $q$ – analogue of (10)

$$\sum_{k=0}^{n} A(2n,2k,q^4,q,q^2)\prod_{j=0}^{k-1} T(2j,q^4,q,q^2)$$

$$=\sum_{k=0}^{n}\frac{[2k+1]}{[n+k+1]}\begin{bmatrix} 2n\\ n-k\end{bmatrix}\frac{1+q^{2k+1}}{1+q^{n+k+1}}\frac{q^{k^2-k}}{(-q;q)_{n-k}(-q;q)_{n+k}}$$

$$=\frac{2}{1+q^{2n}}\begin{bmatrix} 2n\\ n\end{bmatrix}\frac{1}{\prod_{j=1}^{n}(1+q^j)^2}=\frac{2}{1+q^{2n}}A(2n,0,q^2,q,q^2).$$

Another special case of (49) is

$$\sum_{k=0}^{n} A(2n,2k,q^2,q,q^2)\prod_{j=0}^{k-1} T(2j,q^2,q,q^2)=\frac{(-1;q^2)_n}{(-q;q^2)_n}=\prod_{j=0}^{n-1}\frac{1+q^{2j}}{1+q^{2j+1}}.$$